\documentclass[12pt]{amsart}
\usepackage{amsmath}
\usepackage{amsxtra}
\usepackage{amscd}
\usepackage{amsthm}
\usepackage{amsfonts}
\usepackage{amssymb}
\usepackage{eucal}

\input prepictex
\input pictex
\input postpictex

\newtheorem{thm}{Theorem}[section]
\newtheorem{lem}{Lemma}[section]
\newtheorem{cor}{Corollary}[section]
\newtheorem{prop}{Proposition}[section]

\theoremstyle{definition}

\theoremstyle{remark}
\newtheorem{rem}{Remark}[section]

\numberwithin{equation}{section}

\begin{document}

\newcommand{\thmref}[1]{Theorem~\ref{#1}}
\newcommand{\secref}[1]{Section~\ref{#1}}
\newcommand{\lemref}[1]{Lemma~\ref{#1}}
\newcommand{\propref}[1]{Proposition~\ref{#1}}
\newcommand{\corref}[1]{Corollary~\ref{#1}}
\newcommand{\remref}[1]{Remark~\ref{#1}}
\newcommand{\eqnref}[1]{(\ref{#1})}
\newcommand{\exref}[1]{Example~\ref{#1}}

\newcommand{\nc}{\newcommand}
\nc{\on}{\operatorname} \nc{\Z}{{\mathbb Z}} \nc{\C}{{\mathbb C}}
\nc{\oo}{{\mf O}} \nc{\R}{{\mathbb R}} \nc{\N}{{\mathbb N}}
\nc{\bib}{\bibitem} \nc{\pa}{\partial} \nc{\F}{{\mf F}}
\nc{\rarr}{\rightarrow} \nc{\larr}{\longrightarrow}
\nc{\al}{\alpha} \nc{\ri}{\rangle} \nc{\lef}{\langle} \nc{\W}{{\mc
W}} \nc{\gam}{\ol{\gamma}} \nc{\Q}{\ol{Q}} \nc{\q}{\widetilde{Q}}
\nc{\la}{\lambda} \nc{\ep}{\epsilon} \nc{\g}{\mf g} \nc{\B}{\mf B}
\nc{\h}{\mf h} \nc{\Hy}{ \widetilde{H} } \nc{\n}{\mf n}
\nc{\A}{{\mf a}} \nc{\G}{{\mf g}} \nc{\HH}{{\mf h}} \nc{\Li}{{\mc
L}} \nc{\La}{\Lambda} \nc{\is}{{\mathbf i}} \nc{\V}{\mf V}
\nc{\bi}{\bibitem} \nc{\NS}{\mf N}
\nc{\dt}{\mathord{\hbox{${\frac{d}{d t}}$}}} \nc{\E}{\mc E}
\nc{\ba}{\tilde{\pa}}
\def\smapdown#1{\big\downarrow\rlap{$\vcenter{\hbox{$\scriptstyle#1$}}$}}

\nc{\mc}{\mathcal} \nc{\mf}{\mathfrak} \nc{\ol}{\fracline}
\nc{\el}{\ell} \nc{\etabf}{{\bf \eta}} \nc{\x}{{\bf x}}
\nc{\xibf}{{\bf \xi}} \nc{\y}{{\bf y}}

\nc{\NP}{\Pi}

\title{Remarks on the Schur--Howe--Sergeev Duality}

\author[Shun-Jen Cheng]{Shun-Jen Cheng$^\dagger$}
\thanks{$^\dagger$partially supported by NSC-grant 90-2115-M-006-003 of the R.O.C}
\address{Department of Mathematics, National Cheng Kung University, Tainan,
Taiwan}
\email{chengsj@mail.ncku.edu.tw}

\author{Weiqiang Wang}
\address{Department of Mathematics, North Carolina State University,
Raleigh, NC 27695-8205, U.S.A.} \email{wqwang@math.ncsu.edu}

\begin{abstract}
We establish a new Howe duality between a pair of two queer Lie
superalgebras $(q(m),q(n))$. This gives a representation theoretic
interpretation of a well-known combinatorial identity for Schur
$Q$-functions. We further establish the equivalence between this
new Howe duality and the Schur--Sergeev duality between $q(n)$ and
a central extension $\Hy_k$ of the hyperoctahedral group $H_k$. We
show that the zero-weight space of a $q(n)$-module with highest
weight $\la$ given by a strict partition of $n$ is an irreducible
module over the finite group $\Hy_n$ parameterized by $\la$. We
also discuss some consequences of this Howe duality.

\vspace{.3cm}

\noindent{\bf Mathematics Subject Classifications (2000)}: 17B67.
\vspace{.3cm}

\noindent{\bf Key words:} Lie superalgebra $q(n)$,
Schur--Howe--Sergeev duality.

\end{abstract}

\maketitle

\section{Introduction}\label{intro}

This is a sequel to \cite{CW}, in which we studied the Howe
duality between two general linear Lie superalgebras and other
closely related multiplicity-free actions of a general linear Lie
superalgebra, which generalize and unify several classical
results, cf. Howe \cite{H1, H2}.

In this Letter, we construct a new Howe duality involving the
queer Lie superalgebra $q(n)$. The queer Lie superalgebra $q(n)$
(cf., e.g., \cite{K, P}) can be regarded as a true super analog of
the general linear Lie algebra. We show that there is a mutually
centralizing action of $q(m)$ and $q(n)$ on the symmetric
algebra\footnote{In this Letter, we will freely suppress the term
{\em super}. So in case when a superspace is involved, the terms
symmetric, commute etc.~mean {\em super}symmetric, {\em
super}commute etc unless otherwise specified.} of $\C^{mn|mn}$. A
multiplicity-free decomposition of the $q(m)\times q(n)$-module
$S(\C^{mn|mn})$ is explicitly obtained. To achieve this, we use a
remarkable duality, due to Sergeev \cite{Se1} between $q(n)$ and a
finite group $\Hy_k$, which generalizes the celebrated Schur
duality. Here $\Hy_k$ is a central extension of the
hyperoctahedral group.

On the other hand, we show the $(q(m), q(n))$ Howe duality can be
used to re-derive Sergeev duality as well. We also show that the
zero-weight space of a $q(n)$-module with highest weight $\la$,
given by a strict partition of $n$, is an irreducible module over
the finite group $\Hy_n$ parameterized by $\la$. All these are
very much analogous to the classical picture in the general linear
Lie algebra case, cf.~\cite{H2}. As is well known, there is no
unique notion of a Weyl (super)group for a Lie superalgebra. Our
results suggest that $\Hy_n$ may be regarded as a Weyl supergroup
for $q(n)$ in an appropriate sense.

It has been known \cite{Se1} that the characters of the
irreducible $q(n)$-modules under consideration in this letter are
essentially the Schur Q-functions $Q_\la$ (cf. \cite{M}). We
remark that the difficult question of finding the character of a
general finite-dimensional irreducible $q(n)$-module has been
solved recently by Penkov and Serganova \cite{PS}. The $(q(m),
q(n))$-duality can now be interpreted as a representation
theoretic realization of the following well known identity for
Schur Q-functions (cf. Macdonald \cite{M}):

\begin{equation*}
\prod_{i,j=1}^{\infty}\frac{1+x_iy_j}{1-x_iy_j}
=\sum_{\la}2^{-l(\la)} Q_\la(x)Q_\la(y),
\end{equation*}
where $x = (x_1, x_2, \ldots)$, $y=(y_1,y_2,\ldots)$ and the
summation is over all strict partitions.

In the case when $n=1$, the $(q(m),q(n))$ Howe duality essentially
tells us that the $k$-th symmetric algebra of the natural
representation of $q(n)$ is irreducible of highest weight $(k, 0,
\ldots, 0)$. When $m=n$, the $(q(m),q(n))$ Howe duality implies
the existence of a distinguished basis for the center of the
universal enveloping algebra of $q(n)$ (also compare \cite{Se2})
parameterized by strict partitions of length not exceeding $n$. It
is a very interesting question to give a more precise description
of this basis and its relation with symmetric functions as $m$ goes
to infinity. The results of \cite{CW} and the present work also
suggest that there are other Howe dual pairs involving various Lie
superalgebras (cf.~\cite{K}) which deserve further study.

The plan of the letter is as follows.  In \secref{prelim} we
recall some representation-theoretic background of the queer Lie
superalgebras with emphasis on Schur--Sergeev  duality.
\secref{qq-dual} is devoted to establishing the
$(q(m),q(n))$-duality and to the study of its consequences.

\section{The Schur--Sergeev duality}\label{prelim}

Let $\C^{m|n}$ denote the complex vector superspace of dimension
$m|n$, and $\frak{gl}(m|n)$ the Lie superalgebra of linear
transformations of $\C^{m|n}$ (see, e.g.,~\cite{K}). Choosing a
homogeneous basis of $\C^{m|n}$ we may regard $\frak{gl}(m|n)$ as
the space of complex $(m+n)\times(m+n)$ matrices. In the case when
$m=n$ consider an odd automorphism $P:\C^{m|m}\rightarrow
\C^{m|m}$ with $P^2=-1$. The linear transformations of
$\frak{gl}(m|m)$ preserving $P$ is a subalgebra of
$\frak{gl}(m|m)$, denoted by $q(m)$. We have
$q(m)=q(m)_{\bar{0}}\oplus q(m)_{\bar{1}}$, with $q(m)_{\bar{0}}$
isomorphic to the general linear Lie algebra $\frak{gl}(m)$ and
$q(m)_{\bar{1}}$ isomorphic to the adjoint module of
$\frak{gl}(m)$ (cf.~\cite{K}, \cite{P}). Choosing $P$ to be the
$2m\times 2m$ matrix
\begin{equation} \label{oddauto}
\begin{pmatrix}0&I\\ -I&0\end{pmatrix}
\end{equation}
with $I$ denoting the identity $m\times m$ matrix, we may identify
$q(m)$ inside $\frak{gl}(m|m)$ with the space of complex $2m\times
2m$ matrices of the form:
\begin{equation}\label{q(m)}
\begin{pmatrix}
A&B\\
B&A\\
\end{pmatrix},
\end{equation}
where $A$ and $B$ are arbitrary complex $m\times m$ matrices.  Of
course the even elements of $q(m)$ are those for which $B=0$,
while the odd elements are those for which $A=0$.

We will below recall some aspects of finite-dimensional
irreducible representations of $q(m)$ (cf.~\cite{P}).  Let $\B$ be
the Borel subalgebra consisting of those matrices in \eqnref{q(m)}
with $A$ and $B$ upper triangular. Furthermore, let $\NS$ be the
nilpotent subalgebra of $\B$ consisting of those matrices in
\eqnref{q(m)} with $A$ and $B$ strictly upper triangular.  The
Cartan subalgebra $\HH$ is the subalgebra of $q(m)$ consisting of
those matrices in \eqnref{q(m)} with $A$ and $B$ diagonal.  Taking
a linear form $\la$ on $\HH_{\bar{0}}$ we may consider the
symmetric bilinear form on $\HH_{\bar{1}}$ defined by
$(a|b)_{\la}:=\la([a,b])$, $a,b\in\HH_{\bar{1}}$.  Now if
$\HH'_{\bar{1}}\subset\HH_{\bar{1}}$ is a maximal isotropic
subspace with respect to this bilinear form, we may extend $\la$
to a one-dimensional representation of
$\HH_{\bar{0}}+\HH'_{\bar{1}}$ in a trivial way. Inducing from
this we obtain an irreducible $\HH$-module.  This module has an
odd automorphism if and only if the dimension of the quotient
space $\HH_{\bar{1}}/{\rm ker}(\cdot|\cdot)_{\la}$ is odd.  We now
can extend this irreducible $\HH$-module to an irreducible
$\B$-module by letting $\NS$ act trivially. This way one obtains
all finite-dimensional irreducible $\B$-modules. Inducing further
we obtain the Verma module of $q(m)$ associated to the linear form
$\la\in\HH^*_{\bar{0}}$ and, thus, an irreducible $q(m)$-module of
highest weight $\la$ by dividing by its maximal proper submodule.
We note that an odd automorphism of an irreducible $\HH$-module
descends to an odd automorphism of this irreducible quotient. Thus
one may associate to an $m$-tuple $\la=(a_1,\ldots,a_m)$ of
complex numbers an irreducible representation $U_m^{\la}$ of
highest weight $\la$, which is finite-dimensional if and only if
$a_i-a_{i+1}\in\Z_+$ and $a_i=a_{i+1}$ implies that $a_i=0$ for
all $i=1,\ldots,m-1$, cf. \cite{P}.

Recall a partition $\la =(\la_1, \la_2, \ldots, \la_l)$ of length
$l$ is called {\em strict} if $\la_1>\la_2>\cdots > \la_l >0.$ We
will identify $\la $ with $(\la_1, \la_2, \ldots, \la_l, 0,
\ldots, 0 )$ by adding zeros in the end and denote by $|\la|$ the
sum $\la_1 +\cdots +\la_l$. We see that a partition $\la$ may be
regarded as a highest weight of a finite-dimensional irreducible
$q(m)$-module if and only if $\la$ is {\em strict} with length
$l(\la)$ not exceeding $ m$. Furthermore, we have \cite{P}:
\begin{equation}\label{octahedral1}
{\rm dim}\;({\rm
Hom}_{q(m)}(U_m^\la,U_m^\mu))=\delta_{\la\mu}2^{\delta(l(\la))},
\end{equation}
where ${\rm Hom}$ is to be understood in the $\Z_2$-graded sense,
and the number $\delta(l(\la))$ is $0$ for $l(\la)$ even and $1$
otherwise.  More precisely, ${\rm
Hom}_{q(m)}(U_m^{\la},{U}_m^{\mu}))$ is isomorphic to $\C$ in the
case when $l(\la)$ is even and it is isomorphic to a Clifford
superalgebra in one odd variable in the case $l(\la)$ is odd. We
also recall that the character ${\rm ch}U^\la_{m}$ is defined as
the trace of the diagonal matrix ${\rm
diag}(x_1,\ldots,x_m;x_1,\ldots,x_m)$ in $q(m)$ acting on
${U^\la_{m}}$.

We denote by $\NP_k$ the group generated by
$1,z,a_1,\ldots,a_k$ subject to the relations
\begin{equation*}
z^2=1,\ a_i^2=z \ {\rm and}\  a_ia_j= z a_ja_i,\ {\rm for}\
i\not=j.
\end{equation*}
The symmetric group $S_k$ acts on $\NP_k$ via permutations of the
elements $a_1,\ldots,a_k$ and we may thus form the semidirect
product $S_k\ltimes\NP_k$. We denote this semidirect product
(again a finite group) by $\Hy_k$ which is naturally $\Z_2$-graded
by putting $p(a_i)=\bar{1}$, $p(z)=\bar{0}$, and
$p(\sigma)=\bar{0}$ for $\sigma\in S_k$.   Thus $\Hy_k$ is a
finite supergroup in the sense of \cite{J}, i.e.~this is a group
with a subgroup of index two, whose elements we call even and by
definition all other elements are odd (see \remref{supergroup}).
In this letter we will only concern about the {\em $\Z_2$-graded
spin modules\footnote{In \cite{J, J2} these are called {\em
negative supermodules}.}} of $\Hy_k$ (i.e. those $\Z_2$-graded
modules on which $z$ acts as $-1$) which are equivalently modules
over the group superalgebra $\B_k =\C[\Hy_k]/\langle z=-1
\rangle$.

\begin{rem}\label{supergroup}
For a supergroup $G$ by $G$-module homomorphisms we will mean
linear maps commuting with the respective group action. Explicitly
this means that for a homogeneous linear map $f:V\rightarrow W$,
where $V$ and $W$ are supermodules, to be a $G$-homomorphism we
must have $f(g\cdot v)=(-1)^{p(f)p(g)}g\cdot f(v)$, for all $v\in
V$. Recall \cite{J} that given two modules $V$ and $W$ over a
supergroup $G$, the space ${\rm Hom}_{\C}(V,W)$ is naturally a
$G$-module via the action $(g\cdot
T)(v):=(-1)^{p(g)p(T)}gT(g^{-1}v)$, which amounts to giving
$V\otimes W$ a $G$-module structure via the action
$g\cdot(v\otimes w)=(-1)^{p(g)(p(v)+p(w))}gv\otimes gw$. Thus in
what follows the action on a tensor product of two $G$-modules
will always be given this $G$-module structure. By $G$-invariants
inside $V\otimes W$ we shall always mean the usual invariants,
i.e.~$(V\otimes W)^G=\{\gamma\in V\otimes
W|g\cdot\gamma=\gamma,\forall g\in G\}$.
\end{rem}

According to \cite{Se1} and \cite{J}, the ($Z_2$-graded)
irreducible spin modules of $\Hy_k$ are also parameterized
by strict partitions. For strict partitions $\la$ and $\mu$ let
$T_k^\la$ and $T_k^\mu$ denote the corresponding irreducible spin
modules over $\Hy_k$.  We have (\cite{J}, \cite{Se1}) :
\begin{equation}\label{octahedral2}
{\rm dim}({\rm
Hom}_{\Hy_k}(T_k^\la,T_k^\mu))=\delta_{\la\mu}2^{\delta(l(\la))},
\end{equation}
where ${\rm Hom}$ is again to be understood in the $\Z_2$-graded
sense. Furthermore it is known (cf. \cite{J2}) that the character
value of $T_k^\la$ is real and thus $T_k^\la$ is
self-contragredient.

Let us now consider the natural action of $q(m)$ on $\C^{m|m}$. We
may form the $k$-fold tensor product ${\bigotimes}^k\C^{m|m}$, on
which $q(m)$ acts naturally.  In addition we have an action of the
finite supergroup $\Hy_k$: the symmetric group in $k$ letters acts
on ${\bigotimes}^k\C^{m|m}$ by permutations of the tensor factors
with appropriate signs (corresponding to the permutations of odd
elements in $\C^{m|m}$). However, we also have an action of $a_i$
on ${\bigotimes}^k\C^{m|m}$ by means of exchanging the parity of
$i$-th copy of $\C^{m|m}$ via the odd automorphism of $\C^{m|m}$
given by the matrix $P$ of \eqnref{oddauto}.
More explicitly, $a_i$ transforms the vector $v_1\otimes\ldots
v_{i-1}\otimes v_i\otimes\ldots\otimes v_k$ in
$\bigotimes^k\C^{m|m}$ into
$(-1)^{p(v_1)+\ldots+p(v_{i-1})}v_1\otimes\ldots v_{i-1}\otimes
P(v_i)\otimes\ldots\otimes v_k$.

The following remarkable theorem is due to Sergeev \cite{Se1},
which will be referred to as (Schur--)Sergeev duality throughout
the Letter. We refer the reader to \cite{M} for definitions and
properties of the Schur $Q$-functions $Q_{\la}$.

\begin{thm} [Sergeev] \label{Sergeev2} The actions of $q(m)$ and
$\Hy_k$ on the space $\bigotimes^k C^{m|m}$ commute and
$\bigotimes^k\C^{m|m}$ is completely reducible over
$q(m)\times\Hy_k$. Explicitly we have
\begin{equation*}
\bigotimes^k \C^{m|m} \cong \sum_{\lambda}
2^{-\delta(l(\la))}U^{\lambda}_{m} \otimes T_k^{\lambda},
\end{equation*}
where $\lambda$ is summed over all strict partitions with
$|\la|=k$ and $l(\la)\le m$. Furthermore, the character ${\rm
ch}U_m^{\la}$ is given by
$2^{\frac{\delta(l(\la))-l(\la)}{2}}Q_\la(x)$.
\end{thm}

\begin{rem}
The expression $2^{-\delta(l(\la))}U^{\lambda}_{m} \otimes
T_k^{\lambda}$ above has the following meaning. Suppose $A$ and
$B$ are two superalgebras and $V_A$ and $V_B$ are irreducible
modules over $A$ and $B$ such that ${\rm Hom}_{A}(V_A,V_A)$ and
${\rm Hom}_B(V_B,V_B)$ are both isomorphic to the Clifford
superalgebra in one odd variable.  We remark that in the language
of \cite{Se1} (respectively \cite{J}) this is to say that $V_A$
and $V_B$ are irreducible, but not absolutely irreducible
(respectively are of $Q$-type). It is known that $V_A\otimes V_B$
as a module over $A\otimes B$ is not irreducible, but decomposes
into a direct sum of two isomorphic copies (via an odd
isomorphism) of the same irreducible representation (see, e.g.,
\cite{C, J}). In our particular setting when $l(\la)$ is odd both
$T_k^\la$ and $U_m^\la$ are such modules by \eqnref{octahedral1}
and \eqnref{octahedral2}. So in this case we mean to take one copy
inside their tensor product.
\end{rem}
\section{The $(q(m),q(n))$ Howe duality}\label{qq-dual}
Recall that $\NP_k$ acts on $\bigotimes^k\C^{m|m}$ and
$\bigotimes^k\C^{n|n}$ and hence the diagonal subgroup $\Delta
\NP_k\subset \NP_k\times \NP_k$ acts on their tensor product
$$(\bigotimes^k\C^{m|m})\otimes(\bigotimes^k\C^{n|n})
\cong\bigotimes^k(\C^{m|m}\otimes\C^{n|n})$$
(see \remref{supergroup}). So does the symmetric group $S_k$. This
gives rise to the diagonal action of $\Hy_k$.

\begin{lem} \label{generic6}
As a module over $q(m)\times q(n)$, we have
$$(\bigotimes^k(\C^{m|m}\otimes\C^{n|n}))^{\Delta (\Hy_k)}\cong
S^k(\C^{mn|mn}),$$ where $(\cdot)^{\Delta (\Hy_k)}$ denotes the
space ${\Delta (\Hy_k)}$-invariants.
\end{lem}

\begin{proof}  Recall
that $P$ is the odd automorphism given by the matrix
\eqnref{oddauto}. Consider first
$(\C^{m|m}\otimes\C^{n|n})^{\Delta P}$, that is, the $\Delta(P) (=
P \times P)$-invariants in $\C^{m|m}\otimes\C^{n|n}$. Letting
$v_m\in\C^{m|0}$ and $v_n\in\C^{n|0}$, it is clear that this space
consists of elements of the form $v_m\otimes v_n + P(v_m)\otimes
P(v_n)$ and $v_m\otimes P(v_n)+P(v_m)\otimes v_n$, and hence is
isomorphic to $\C^{mn|mn}$. Therefore, since $\Delta \NP_k$ is a
subgroup (since all elements are now even) of $\Delta(S_k\ltimes
\NP_k)$ generated by the $k$ copies of $\Delta(P)$'s, we have
\begin{align*}
& (\bigotimes^k\C^{m|m}\otimes\C^{n|n})^{\Delta (S_k\ltimes
\NP_k)}\\
 & \cong ((\bigotimes^k\C^{m|m}\otimes\C^{n|n})^{\Delta
\NP_k})^{\Delta S_k} \\ & \cong(\bigotimes^k(
(\C^{m|m}\otimes\C^{n|n})^{\Delta P}))^{\Delta S_k} \\ &
\cong(\bigotimes^k\C^{mn|mn})^{S_k} \\ & \cong S^k(\C^{mn|mn}).
\end{align*}
\end{proof}

Let $x_1^i,\ldots,x_m^i$, $\xi_1^i,\ldots,\xi_m^i$, for
$i=1,\ldots,n$ denote the standard coordinates of $\C^{mn|mn}$. We
may then identify $S(\C^{mn|mn})$ with the polynomial algebra
generated by $x_i$ and $\xi_j$.  Introduce the following first
order differential operators:
\begin{align}
A^{pq}&=\sum_{i=1}^m(x_i^p\frac{\partial}{\partial
x_i^q}+\xi_i^p\frac{\partial}{\partial \xi_i^q}),\quad 1\le p,q\le
n,\nonumber\\
B^{pq}&=\sum_{i=1}^m(x_i^p\frac{\partial}{\partial
\xi_i^q}-\xi_i^p\frac{\partial}{\partial x_i^q}),\quad 1\le p,q\le
n,\label{generic3}\\
A_{pq}&=\sum_{j=1}^n(x^j_p\frac{\partial}{\partial
x^j_q}+\xi^j_p\frac{\partial}{\partial \xi^j_q}),\quad 1\le p,q\le
m,\nonumber\\
B_{pq}&=\sum_{j=1}^n(x^j_p\frac{\partial}{\partial
\xi^j_q}+\xi^j_p\frac{\partial}{\partial x^j_q}),\quad 1\le p,q\le
m.\label{generic4}
\end{align}

The following lemma can be proved directly.

\begin{lem}
The operators $A_{pq}$ and $B_{pq}$, for $1\le p,q\le m$, form a
copy of $q(m)$, while $A^{pq}$ and $B^{pq}$, for $1\le p,q\le n$,
form a copy of $q(n)$. Furthermore, they define a commuting action
of $q(m)$ and $q(n)$ in $S(\C^{mn|mn})$.
\end{lem}

\begin{thm}\label{qq-duality}
The action of $q(m)\times q(n)$ on $S(\C^{mn|mn})$ is
multiplicity-free.  More precisely we have the following
decomposition:
\begin{equation*}
S^k(\C^{mn|mn})\cong \sum_{\la}2^{-\delta(l(\la))} U_m^\la\otimes
U_n^\la,
\end{equation*}
where $\la$ is summed over all strict partitions of length not
exceeding ${\rm min}(m,n)$.
\end{thm}

\begin{proof}
By Schur--Sergeev duality (\thmref{Sergeev2}) we have
$${\bigotimes^k}{\C^{m|m}}\cong\sum_{\la}2^{-\delta(l(\la))}U^\la_m\otimes
T^\la_k,$$ as $q(m) \times \Hy_k$ module, where the summation is
over strict partitions $\la$ with length $l(\la)\le m$.  Therefore
combined with \lemref{generic6} this gives us
\begin{align*}
S^k(\C^{mn|mn})&
\cong((\bigotimes^k\C^{m|m})\otimes(\bigotimes^k\C^{n|n}))^{\Delta
(S_k\ltimes \NP_k)}\\
&\cong
\sum_{\la,\mu}2^{-\delta(l(\la))}2^{-\delta(l(\mu))}(U_m^\la\otimes
T_k^\la\otimes U_n^\mu\otimes T_k^\mu)^{\Delta (S_k\ltimes
\NP_k)}\\ &\cong
\sum_{\la,\mu}2^{-\delta(l(\la))}2^{-\delta(l(\mu))}(U_m^\la\otimes
U_n^\mu)\otimes (T_k^\la\otimes T_k^\mu)^{\Delta (S_k\ltimes
\NP_k)}.
\end{align*}
Now by \eqnref{octahedral2} and the fact that irreducible $\Hy_k
(= S_k\ltimes \NP_k)$-modules are self-contragredient we have
\begin{align*}
S^k(\C^{mn|mn})&\cong
\sum_{\la}2^{\delta(l(\la))}2^{-2\delta(l(\la))}U_m^\la\otimes{U}_m^\la\\
&\cong \sum_{\la}2^{-\delta(l(\la))}U_m^\la\otimes U_m^\la,
\end{align*}
where the summation is over all strict partitions of length $\le {\rm
min}(m,n)$.
\end{proof}

Comparing the characters of both sides of the
$(q(m),q(n))$-duality (\thmref{qq-duality}) we obtain, by using
\thmref{Sergeev2}, that
\begin{equation*}
\prod_{i=1}^m\prod_{j=1}^{n}\frac{1+x_iy_j}{1-x_iy_j}
=\sum_{\la}2^{-\delta(l(\la))}
2^{\frac{\delta(l(\la))-l(\la)}{2}}Q_\la(x)
2^{\frac{\delta(l(\la))-l(\la)}{2}}Q_\la(y),
\end{equation*}
which is equivalent to
\begin{equation*}
\prod_{i=1}^m\prod_{j=1}^{n}\frac{1+x_iy_j}{1-x_iy_j}
=\sum_{\la}2^{-l(\la)} Q_\la(x)Q_\la(y),
\end{equation*}
where the summation is over all strict partitions of length $\le
{\rm min}(m,n)$.  This identity of Schur $Q$-functions is well
known (see e.g.~\cite{M}) and \thmref{qq-duality} provides a
representation-theoretic interpretation of it.

The next corollary is immediate from \thmref{qq-duality}.

\begin{cor}
The images of the universal enveloping algebras of $q(m)$ and
$q(n)$ in the endomorphism algebra of $S^k(\C^{mn|mn})$ are mutual
centralizers.
\end{cor}

When $n=1$ the $(q(m),q(n))$-duality reads
\begin{equation*}
S^k(\C^{m|m})\cong\frac{1}{2} (U_m^{(k)}\otimes U_1^{(k)}),
\end{equation*}
where $(k)$ above denotes the one-part partition. Since
$U_1^{(k)}$ is a two-dimensional module, the right-hand side is
exactly $U_m^{(k)}$ as a $q(m)$-module. Hence we have established
the following.

\begin{prop}\label{generic5}
The $k$-th symmetric tensor of $\C^{m|m}$ is the irreducible
$q(m)$-module $U_m^{(k)}$, associated to the one-part partition
$(k)$.
\end{prop}

\begin{rem} Another proof of \propref{generic5} goes as follows.
Recalling that $S^k(\C^{m|n})=\bigoplus_{i=0}^k
S^i\C^m\otimes\Lambda^{k-i}\C^n$, we have therefore ${\rm
ch}S^{k}(\C^{m|m})=\sum_{i=0}^k h_ie_{k-i}=q_k$ (see \cite{M}
pp.~261 for notation), which coincides with $Q_{(k)}$. Here $h_i$
and $ e_i$ denote the $i$-th complete and respectively elementary
symmetric functions. But ${\rm ch}U^{(k)}_m=Q_{(k)}$ by \cite{Se1}
and so $S^{k}(\C^{m|m})\cong U^{(k)}_m$.
\end{rem}

\begin{rem}\label{generic1}
Given a left module $M$ over a Lie superalgebra $\G$, we can make
$M$ into a right module by defining $m\cdot
x:=-(-1)^{p(m)p(x)}xm$, for $m\in M$ and $x\in\G$.  Now if $M$ in
addition has a left module structure over another Lie superalgebra
$\G'$ such that the action of $\G$ and $\G'$ commute, then the so
induced right action of $\G$ on M will not commute with the left
action of $\G$ (in the super-sense), but rather they will commute
with each other in the usual sense. Thus the induced right action
of $q(n)$ above has the result that it commutes with the left
action of $q(m)$ in the usual sense.
\end{rem}

Let $GL(n)$ be the Lie group whose Lie algebra is the even part of
the Lie superalgebra $q(n)$ and let $A_n$ denote its diagonal
torus. Given a $q(n)$-module (or a $GL(n)$-module) $U$, we call
the subspace $U^{A_n, {\rm det}}$ inside $U$, which transforms
under the action of $A_n$ by the determinant character, the {\em
zero-weight space} of $U$.

\begin{thm} \label{th_equiv}
The $(q(m),q(n))$ Howe duality implies the Sergeev duality.
\end{thm}

\begin{thm} \label{th_zerowt}
Given a strict partition $\la$ of $n$, the zero weight space of
$U^{\la}_n$ admits a natural action of the finite group $\Hy_n$,
and it is isomorphic to the irreducible module $T^{\la}_n$.
\end{thm}

\begin{proof} We will establish \thmref{th_equiv}
and \thmref{th_zerowt} together. The argument we will present
follows closely the one used in \cite{H2} to derive the Schur
duality from the $(\frak{gl}(m),\frak{gl}(n))$-duality.

The $(q(m),q(n))$-duality says that $S(\C^{m|m}\otimes\C^n)\cong
\sum_{\la}2^{-\delta(l(\la))}{U}^\la_m\otimes{U}_n^\la$, where the
summation is over all strict $\la$ with $l(\la)\le{\rm min}(m,n)$.
Observe that the space
$\bigotimes^n\C^{m|m}\cong\C^{m|m}\otimes\C^n$ may be identified
with the zero-weight space $S(C^{m|m}\otimes\C^n)^{A_n,{\rm
det}}$. Putting these together we have
\begin{equation} \label{eq_side}
\bigotimes^n\C^{m|m}\cong\sum_{\la}2^{-\delta(l(\la))}
{U}_m^\la\otimes({U}_n^\la)^{A_n,{\rm det}},
\end{equation}
where $\la$ runs over all strict partitions $n$. So to recover
\thmref{Sergeev2} it suffices to show that $({U}_n^\la)^{A_n,{\rm
det}}$ is isomorphic to the irreducible $\Hy_n$-module $T^{\la}_n$,
which is the contents of \thmref{th_zerowt}.

First note that the normalizer of $A_n$ acts on the $A_n$-weight
spaces by permutation.  Since the determinant character is
invariant under permutation of weights, we see that
$({U}_n^\la)^{A_n,{\rm det}}$ is invariant under the action of
$S_n$.  Now the induced right action (see \remref{generic1}) of
the operators $B^{ii}$ (see \eqnref{generic3}) when acting on
$({U}_n^\la)^{A_n,{\rm det}}$ satisfy the commutation relations of
the $a_i$'s in $\Hy_n$, for $i=1,\ldots,n$.  This action of the
$B^{ii}$'s combined with the action of $S_n$ then gives a right
action of $\Hy_n$ on $({U}_n^\la)^{A_n,{\rm det}}$.

Set $n=m$ in the remainder of the proof. (However, it is
convenient to continue making the distinction between $n$ and
$m$.) Consider
\begin{equation} \label{eq_reg}
(\bigotimes^n\C^{m|m})^{A_m,{\rm det}}\cong
\sum_{\la}2^{-\delta(l(\la))}({U}_m^\la)^{A_m,{\rm
det}}\otimes({U}_n^\la)^{A_n,{\rm det}}.
\end{equation}
Recall that $x_1^i,\ldots,x_m^i$, $\xi_1^i,\ldots,\xi_m^i$
$(i=1,\ldots,n)$ are the standard coordinates of $\C^{mn|mn}
\cong\C^{m|m}\otimes\C^n$. It is not difficult to see that the
space $(\bigotimes^n\C^{m|m})^{A_m,{\rm det}}$ inside
$S(\C^{m|m}\otimes\C^n)^{A_m,{\rm det}}$, which is
$S(\C^{m|m}\otimes\C^n)^{A_m \times A_n,{\rm det} \times {\rm
det}}$, may be identified with the space spanned by vectors of the
form $v_1^{\sigma_1}v_2^{\sigma_2}\ldots v_n^{\sigma_n}$, where
$v$ denotes either $x$ or $\xi$ and $\sigma$ is a permutation of
$\{1,\ldots,n\}$.  Thus this space is in bijection with the space
$\B_n=\C[\Hy_n]/\langle z=-1 \rangle$.  On this space the
normalizer of $A_m$ acts, and so we obtain an action of the
symmetric group $S_m$.  We also have an action of $B_{jj}$ (see
\eqnref{generic4}), which gives rise to the action of $a_j$ in
$\Hy_m$, for $j=1,\ldots,m$.  This action combined with that of
$S_m$ gives a left action of $\Hy_m$ on $(U_m^\la)^{A_m,\rm det}$.
Furthermore the induced right action of $\Hy_n$ above and this
action commute in the usual sense according to \remref{generic1}.

Also our left action of $S_n$ permutes the upper indices of the
vector $v_1^{\sigma_1}v_2^{\sigma_2}\ldots v_n^{\sigma_n}$,
whereas the left action of the $a_i$'s changes the parity of
$v_i^{\sigma_i}$. Thus this is the left regular action of $\Hy_n$.
Our right action of $S_n$, on the other hand, permutes the lower
indices of $v_1^{\sigma_1}v_2^{\sigma_2}\ldots v_n^{\sigma_n}$,
whereas our right action of $a_i$ changes the parity of
$v_i^{\sigma^{-1}_i}$. Thus our right action is the right regular
representation of $\Hy_n$. From the general theory of finite
supergroup \cite{J}, the left-hand side of \eqnref{eq_reg}, which
is isomorphic to $\B_n$ under the left and right actions of
$\Hy_n$, is equal to the summation
$\sum_{\la}2^{-\delta(l(\la))}{T}^\la_n \otimes{T}^\la_n$ over all
strict partitions of $n$. Decomposing $({U}_m^\la)^{A_m,{\rm
det}}$ in the right-hand side of \eqnref{eq_reg} into a direct sum
of irreducible $\Hy_m$-modules, we see that this is only possible
when each $({U}_m^\la)^{A_m,{\rm det}}$ itself is irreducible as a
$\Hy_m$-module. Comparing with \eqnref{eq_side}, we see that this
irreducible module is isomorphic to $T^{\la}_n$.
\end{proof}

In \cite{Se2} Sergeev computed the center of the universal
enveloping algebra $U(q(m))$ of the Lie superalgebra $q(m)$.  We will see
that a different description of it can also be obtained from the
$(q(m),q(n))$-duality.

Recall that for a finite-dimensional Lie superalgebra $\G$ we have
$S(\G)^{\G}\cong Z(U(\G))$ as a $\G$-module (cf. \cite{Se2}), the
center of the universal enveloping algebra of $\G$.  Now replacing
$q(n)$ and $\C^{n|n}$ in the $(q(m),q(n))$-duality by $q(m)$ and
$\C^{m|m*}$, the action contragredient to the natural action of
$q(m)$, we have
\begin{equation*}
S^k(\C^{m|m}\otimes\C^{m|m*})^{\NP_k}
\cong2^{-\delta(l(\la))}\sum_{\la}{U}_m^\la\otimes{U}^{\la *}_m.
\end{equation*}
However, from our earlier description of
$(\C^{m|m}\otimes\C^{m|m*})^{\Delta P}$, we see that, as a
$q(m)$-module, it is isomorphic to the adjoint representation of
$q(m)$.  Thus
\begin{align*}
S(q(m))^{q(m)}&\cong
\sum_{\la}2^{-\delta(l(\la))}({U}_m^\la\otimes{U}_m^{\la
*})^{\Delta q(m)}\\ &\cong\sum_{\la}2^{-\delta(l(\la))}{\rm
Hom}_{q(m)}({U}_m^\la,{U}_m^{\la }),
\end{align*}
where the summation is over all strict partitions $\la$ with
$l(\la)\le m$.  Thus combined with \eqnref{octahedral1} we have
proved the following proposition.

\begin{prop}
The center of the universal enveloping algebra of $q(m)$ admits a
distinguished basis parameterized by strict partitions of length
less than or equal to $m$.
\end{prop}

{\bf Acknowledgment.} After we completed this work, we came across
a preprint of Sergeev,`An analog of the classical invariant theory
for Lie superlagebras', math.RT/9810113, where he independently
obtained the $(q(m), q(n))$ Howe duality (i.e. \thmref{qq-duality}
in this Letter). The other results of this Letter seem to be new.

\end{document}